\documentclass[a4paper,10pt,reqno]{amsart}        
\usepackage{verbatim}       
\usepackage{amssymb}        
\usepackage{enumerate}       
\usepackage[active]{srcltx}     
\usepackage{float}  
\numberwithin{equation}{section}

\makeatletter
\@namedef{subjclassname@2020}{%
  \textup{2020} Mathematics Subject Classification}
\makeatother

\usepackage{tikz}
\usetikzlibrary{positioning}
\usetikzlibrary{patterns}
\usepackage{marvosym}
\usetikzlibrary{arrows.meta}
\usetikzlibrary{decorations.text}

\usepackage{t1enc}       
\usepackage[utf8x]{inputenc}       
\usepackage{tikz}

\usetikzlibrary{shapes.multipart}

\newif\ifdeveloping

\ifdeveloping
\usepackage[inline]{showlabels}
\showlabels[\small\color{gray}]{cite}
\showlabels[\small\color{red}]{bibitem}
\fi

\newtheorem{theorem}{Theorem}[section]        
\newtheorem{lemma}[theorem]{Lemma} 
\newtheorem{klemma}[theorem]{Key Lemma}       
\newtheorem*{theorem*}{Theorem}

\newtheorem{proposition}[theorem]{Proposition}

\newtheorem{problem}[theorem]{Problem}

\newtheorem{corollary}[theorem]{Corollary}       
\newtheorem{claim}{Claim}[theorem]

\theoremstyle{definition}       
\newtheorem{definition}[theorem]{Definition}       
\theoremstyle{remark}

\newcommand{\mc}[1]{\mathcal{#1}}       
\newcommand{\mbb}[1]{\mathbb{#1}}       
       
\newcommand{\mf}[1]{\mathfrak{#1}}       
\newcommand{\setm}{\setminus}       
\newcommand{\empt}{\emptyset}       
\newcommand{\subs}{\subset}       
\newcommand{\dom}{\operatorname{dom}}       
\newcommand{\ran}{\operatorname{ran}}

\newcommand{\sss}{\operatorname{s}}

\newcommand{\Rhered}{R-hereditarily\ }

\newcommand{\accu}[2]{\operatorname{acc}(#1,#2)}

\def\<{\left\langle}       
\def\>{\right\rangle}

\newcommand{\CD}{\operatorname{CD}}

\newcommand{\PNC}{P-NC}

\newcommand{\DRC}{DRC}

\newcommand{\DCC}{DCC}
\newcommand{\CC}{CC}

\newcommand{\AC}{AC}
\newcommand{\SP}{SP}

 \newcommand{\bluetwenty}{black!0}

 \newcommand{\roomwidth}{4 cm}
 \newcommand{\ctimewidth}{4 cm}
 
 \def\tiii{
  \tikzstyle{zfc}=[draw, rectangle,  minimum height=1cm, minimum width=\roomwidth,anchor=south west]

  \tikzstyle{feltetel}=[draw, rectangle, minimum height=1 cm, minimum width=\ctimewidth,anchor=north east]
  
 } 
 
 \newcommand{\zfckek}[4]
 {\node[zfc,minimum height=#3 cm,fill=\bluetwenty]  at (#1,#2)   {{#4}}; 
 }

 \newcommand{\zfcpiros}[4]
 {\node[zfc,minimum height=#3 cm,fill=\bluetwenty]  at (#1,#2)   {{#4}}; 
 }

 \newcommand{\conkek}[4]
 {\node[zfc,minimum height=#3 cm, pattern = dots, pattern color=gray!50, ]  at (#1,#2)   {{#4}}; 
 }

 \newcommand{\ujcon}[4]
 {\node[zfc,minimum height=#3 cm,draw=black, pattern color=gray!50, pattern = north east lines ]  at (#1,#2)   {{#4}}; 
 }

 \newcommand{\ujzfc}[4]
 {\node[zfc,minimum height=#3 cm,draw=black, fill=gray!30  ]  at (#1,#2)   {{#4}}; 
 }

\author{Istv{\'a}n Juh\'asz}
\address{HUN-REN Alfr\'ed Rényi Institute of Mathematics}
\email{juhasz@renyi.hu}
 \author[L. Soukup]{Lajos Soukup}       
 \address       
       { HUn-REN Alfr{\'e}d R{\'e}nyi Institute of Mathematics       
 }       
 \email{soukup@renyi.hu}       
       
 \author[Z. Szentmikl\'ossy]{Zolt\'an Szentmikl\'ossy}       
 \address{Eötvös University of Budapest}       
 \email{szentmiklossyz@gmail.com}       
       
       
       
 \title[Pseudocompact vs  countably compact]{Pseudocompact versus  countably compact in first countable spaces}       
\thanks{The  preparation of this paper was
 supported by NKHT grant   K129211}       
 \date{\today}

 \dedicatory{Dedicated to the memory of Peter Nyikos}
       
\begin{document} 

\begin{abstract} 

  The primary objective of this work is to construct spaces that are ``{\em pseudocompact but not countably compact} '', abbreviated as \PNC, 
  while endowing them with additional properties.

  First, motivated by an old problem of van Douwen concerning first
  countable \PNC\ spaces with countable extent,
  we construct from CH  a locally compact and locally countable first
  countable \PNC\ space with countable spread.

  A space is deemed {\em densely countably compact}, denoted as \DCC\ for brevity, 
  if it possesses a dense, countably compact subspace. Moreover, a space qualifies as 
  {\em densely relatively countably compact}, abbreviated as \DRC, if it contains a
   dense subset $D$
   such that every infinite subset of $D$ has an accumulation point in $X$.

A countably compact space is \DCC, 
   a \DCC\ space is \DRC, 
   and a \DRC\ space is evidently pseudocompact. 
   The Tychonoff plank is a \DCC\ space but is not countably compact. 
   A $\Psi$-space belongs to the class of \DRC\ spaces but is $\neg$\DCC. 
   Lastly, if $p\in {\omega}^*$ is not a P-point, then $T(p)$, 
   representing the type of $p$ in ${\omega}^*$, constitutes a pseudocompact subspace of 
   ${\omega}^*$ that is $\neg$\DRC.

   When considering a topological property denoted as $Q$, 
   we define a space $X$ as ``{\em \Rhered   $Q$}'' 
   if every regular closed subspace of $X$ also possesses property $Q$. 
  The Tychonoff plank and  the $\Psi$-space  are not R-hereditary examples
  for separating the above-mentioned properties. 
   However, the  aforementioned space $T(p)$ is an  R-hereditary   example, 
   albeit not being first countable. 
   
   In this paper we want to find (first countable) examples
  which separates these properties R-hereditarily. 
We have obtained the following result.

\begin{enumerate}[(1)]
\item There is a \Rhered\ ``\DCC, but not countably compact'' space.
\item \label{two}
If CH holds, then there is a \Rhered\ ``\DRC, but $\neg$\DCC'' space.
\item 
If  $\mathfrak{s}=\mathfrak{c}$,  then there is a first countable, \Rhered\ ``pseudocompact , but $\neg$\DRC'' space.
\end{enumerate}
In contrast to \eqref{two}, it is unknown whether 
a first countable, \Rhered\ ``\DRC, but $\neg$\DCC'' space $X$ can exist.
\end{abstract}       

\maketitle

\section{Introduction}

The concept of {\em pseudo-compactness}  was introduced by Hewitt in \cite{He48}.
In \cite{MaPa55} Mardesic and Papic proposed the notion of {\em feebly compact} spaces, 
and they established  that 
a completely regular space is pseudocompact if and only if  it is feebly compact.  

A countably compact  (abbreviated \CC)  Tychonoff space is necessarily pseudo-compact. 
However, the reverse implication does not hold:
both a  $\Psi$-space  and a  Tychonoff plank serve as  simple examples of 
\emph{pseudocompact, but not countably compact} 
(abbreviated \PNC) spaces.

What   weaker conditions lead to a space being pseudocompact?

A space is pseudocompact if it has a {\em dense, countably compact subspace}, (in short, if the space is \DCC).   
For example, ${\omega}_1\times {\omega}$  is a dense, countably compact subspace of the  
Tychonoff plank. Answering  affirmatively a question of   
Mardesic and Papic, in \cite{Ma71} Marjanovic  showed that  
a $\Psi$-space is  pseudo-compact space which is $\neg$\DCC.

Let us  say that a subspace $D$ of a space $X$ is  
\emph{relatively countably compact} 
iff every infinite subset of $D$
has a limit point in $X$. 
If a topological space  \emph{contains a dense, relatively countably compact subset} 
(it is \DRC, in short), then $X$ is clearly pseudocompact.  
For example,  a  $\Psi$-space   is a \DRC\ space because the isolated points form a dense, 
relatively countably compact subset.    

In \cite{DoSh17} Dorantes-Aldama and Shakhmatov introduced the following concept. 
A topological space $X$ is called {\em selectively pseudocompact}  (abbreviated SP) iff given any family 
$\{U_n:n\in {\omega}\}$ of non-empty open sets, it is possible to choose  points $x_n\in U_n$ such that 
the set $\{x_n:n\in {\omega}\}$ has an accumulation point.
Clearly every \DRC\ space is \SP, and  all the \SP\ spaces are feebly compact.

In \cite[Section 2]{Be81}
Berner constructed a dense subspace of $\Sigma(2^{{\omega}_1})$, referred to as ``Berner's $\Sigma$'',
which is ``\SP, but $\neg$\DRC''.

In \cite[Section 5]{Be81} Berner introduced another example: a 0-dimensional, locally countable,
first countable, ``\SP\ but  
$\neg$\DRC'' space of cardinality $\mf c^+$, which will refer to  as ``Berner's monster''.

Ginsburg and Sacks, \cite{GiSa75}, using a result of Frolik,  proved that if $p\in {\omega}^*$ 
is not a P-point, then 
$T(p)$, the type of $p$ in ${\omega}^*$, is a pseudocompact subspace of ${\omega}^*$.
In \cite{Ku78} Kunen constructed a weak P-point $p$ which is not P-point in ZFC,  
and so the pseudocompact space $T(p)$ mentioned above 
is  an {\em anti-countably compact} space, i.e. no countable subset in it has a limit point.
In \cite{Sa86} Shakhmatov constructed arbitrarily large pseudocompact, anti-countably compact spaces 
in ZFC.

The last two results   addressed the following problem:
{\em To what extent can a pseudocompact space deviate from being countably compact?}

Let us observe that some of the examples mentioned so far possess interesting additional properties.
For a given topological property  $Q$, a space $X$ is defined 
  to be  {\em \Rhered   $Q$} 
   if every regular closed subspace of $X$ also has  property $Q$.
For instance, every pseudocompact (\DCC, \DRC, \SP)  space is \Rhered pseudocompact  
(\DCC, \DRC, \SP, respectively) .
The $\Psi$-space is first countable but not \Rhered ``$\neg$\DCC''. 
Similarly, Berner's monster is  first countable, but not \Rhered ``$\neg$\SP'' as
it is locally compact.

On the other hand,   the space $T(p)$ and the space constructed by Shakhmatov 
are  \Rhered ``pseudocompact, but $\neg$\SP,'' but neither of these spaces is  first countable.
Berner's $\Sigma$ is \Rhered ``SP, but $\neg$\DRC'', though   it, too,  is not first countable.

These observations raise the following question: 
{\em Can we find examples that are both first countable and ``R-hereditary '', while being
as far  from being countably compact as possible, in other words, that contain as many  closed discrete countable sets as possible?
Can you find large ``R-hereditary'' examples, in particular, examples of sizes greater than $2^{\omega}$? 
}

Let us observe that 
if $X$ is an ``example'', e.g. , $X$ is ``\DRC, but $\neg$\DCC'', then the disjoint union of $X$ and a compact space is also an example. Consequently,
it is impossible to establish  a cardinality bound for the sizes of spaces that  are ``\DRC\ but $\neg$\DCC''. 

However, the situation is entirely different when considering 
 ``R-hereditary'' examples, as    the disjoint union of an R-hereditary example and a compact space is not an  R-hereditary  example.

\medskip 

A first countable \DCC\ space is countably compact, so we can not expect 
first countable examples separating \CC\ and \DCC.
The Tychonoff plank is a \DCC, but $\neg$\CC\ space of size ${\omega}_1$ , but it is not 
an R-hereditary example.
However, we can construct arbitrarily large ``R-hereditary''   examples in ZFC.

\begin{theorem}\label{tm:ty+}
  For each cardinal ${\kappa}$,
  there is an \Rhered\ ``\DCC\ but $\neg$\CC'' space $X$ with $|X|= {\kappa}^{\omega}$.
  \end{theorem}
 
  \begin{proof}
    Our space $X$ will be a dense subspace of the compact space 
    ${}^{\omega}({\kappa}^++1)$, namely 
  let \begin{displaymath}
  X=\{f\in {}^{\omega}({\kappa}^++1):|\{n: f(n)={\kappa}^+\}|<{\omega}\}.
  \end{displaymath}
  The subspace $Y={}^{\omega}({\kappa}^+)$ of $X$ is  dense  and  countably compact.
  If ${\varepsilon}$ is an elementary open set in $X$, i.e. 
  $\dom({\varepsilon})\in {[{\omega}]}^{<{\omega}}$ and $\ran({\varepsilon})$
  consists of open subsets of ${\kappa}^+$,
  then define $\{f_n:n\in {\omega}\}$ as follows.
  For each $i\in \dom({\varepsilon})$ pick ${\alpha}_i\in {\varepsilon}(i)$ and let 
  \begin{displaymath}
  {f_n(i)}=\left\{\begin{array}{ll}
  {{\alpha}_i}&\text{if $i\in \dom({\varepsilon})$},\\
  {{\kappa}^+}&\text{if $i\in n\setm \dom({\varepsilon})$},\\
  0&\text{otherwise}.
  \end{array}\right.
  \end{displaymath} 
  Then  $\{f_n:n\in {\omega}\}\subs X\cap[{\varepsilon}]$ is closed discrete
in $X$ because it converges to the  function  
$$\{\<i,{\alpha}_i\>:i\in \dom(\varepsilon)\}\cup \{\<n,{\kappa}^+\>:n\in {\omega}\setm \dom({\varepsilon})\}\in {}^{\omega}({\kappa}^++1)\setm X.$$
So $X$ is \Rhered $\neg$\CC.
\end{proof}

A $\Psi$-space is an example of a  first countable space that is   ``\DRC\ but $\neg$\DCC'', but it is not an R-hereditary example. 
In Section \ref{sc:DRC-notDCC}    we prove  
Theorem \ref{tm:drc-notdcc} which directly  implies the   following result:
\begin{theorem}\label{tm:Thm1.2}
(1)
  If CH holds, then there is an \Rhered ``\DRC, but $\neg$\DCC'' space  $X$ of size ${\omega}_1$.
   (2)  It is consistent that CH holds, $2^{{\omega}_1}$ is as large as you wish, and
  there is an \Rhered ``\DRC, but $\neg$\DCC'' space  $X$ of size $2^{{\omega}_1}$. 

  \end{theorem}
We do not have even a  consistent example of a 
 first countable, \Rhered ``\DRC, but $\neg$\DCC'' space.

\medskip

Berner's $\Sigma$ is \Rhered\ ``\SP\ but $\neg$\DRC'', but its character is ${\omega}_1$.
On the other hand,   Berner's monster is a first countable, \SP\ but $\neg$\DRC\ space, but 
it is not an \Rhered example, as it is locally compact.   
In Section \ref{sc:P-notDRC} we will prove the following result (see Theorems \ref{tm:main} and \ref{tm:geju}).

\begin{theorem}
\label{tm:Tm1.3}  If $\mathfrak s=\mathfrak c$, then there is a first countable,
\Rhered\ ``\SP\ but $\neg$\DRC'' space of size  $\mathfrak c$.    
\end{theorem}

The space $T(p)$ and the example of Shakmatov are anti-countably compact, so  
they are \Rhered  ``pseudocompact, but $\neg\SP$''.
A first countable pseudocompact  space is selectively pseudocompact,
 so we can not expect 
first countable examples separating these properties.

\smallskip

Figure \ref{harmadikabra} provides a summary of our findings. 
The symbol \Lightning\ indicates the non-existence of corresponding spaces, 
while $\surd$ denotes the presence of examples with stronger properties in certain cells. 
Examples are presented with slanted line background  when they represent  consistent constructions.
Question mark indicates   the absence of an example.

     \begin{figure}[ht]\label{harmadikabra}
      \renewcommand{\roomwidth}{3.5 cm}

     \scalebox{0.7}
     {\begin{tikzpicture}[x=\roomwidth, y=-1cm, node distance=0 cm,outer sep = 0pt]
     \tiii
     
     \node[zfc] (A) at (1,1) {$\bf \neg \CC\ \land\ \DCC$};
     \node[zfc] (B) at (2,1) {$\bf \neg \DCC\ \land\ \DRC$   };
     \node[zfc] (C) at (3,1) {$\bf \neg \DRC\ \land\ \SP$  };
     \node[zfc] (C) at (4,1) {$\bf \neg \SP\ \land\ P$  };
     
     \foreach  \x/\y  in {1/--,2/\bf R-hereditary,3/{$\bf M_1$},4/{\bf R-hereditary, $\bf M_1$}}
     \node[feltetel] at  (1,\x)   {\y}  ;
     
     \zfckek{1}{2}{1}{Tychonoff plank}
     \zfckek{2}{2}{1}{$\surd$ }
     \zfckek{3}{2}{1}{$\surd$}
     \zfckek{4}{2}{1}{$\surd$}

     \ujzfc{1}{3}{1}{Thm \ref{tm:ty+}}
     \ujcon{2}{3}{1}{\it Thm \ref{tm:Thm1.2}}
     \zfckek{3}{3}{1}{Berner's  $\Sigma$ }
     \zfckek{4}{3}{1}{$T(p)$, Shakmatov}
     
     \zfckek{1}{4}{1}{\Lightning}
     \zfckek{2}{4}{1}{$\Psi$-space}
     \zfckek{3}{4}{1}{Berner's monster}
     \zfckek{4}{4}{1}{\Lightning}
     
     \zfcpiros{1}{5}{1}{\Lightning}
     \conkek{2}{5}{1}{\bf ??}
     \ujcon{3}{5}{1}{\it Thm  \ref{tm:Tm1.3}}
     \zfckek{4}{5}{1}{\Lightning}
     \end{tikzpicture}
     }
     \caption{Examples separating classes of  pseudocompact spaces}
          \end{figure}
     
\medskip

 The actual starting point of our investigation was  a problem posed by van Douwen.
As we remarked,   both a  $\Psi$-space  and a  Tychonoff plank serve as  simple examples of 
\emph{pseudocompact, but not countably compact}. Notably, a $\Psi$-space is first countable but has uncountable  extent, while  the Tychonoff plank has  countable extent,
 but fails to be first countable.  
So it is a  natural question is  whether there are \PNC\  spaces with small extent and countable character?

Eric Van Douwen  and Peter Nyikos  constructed  two distinct examples of  such spaces,
assuming $\mf b={\omega}_1$ 
(as discussed in  \cite[Notes to Section 13]{HandBook-vD}, where Nyikos provided an example) and 
assuming $\mf b=\mf c$ 
(see \cite[Ex. 13.3]{HandBook-vD}), respectively. In \cite[Question 12.5 and 12.6]{HandBook-vD},
van Douwen posed two related questions: 
the first concerns  the minimum cardinality of a first countable P-NC space, known  
to lie between $\mathfrak{b}$ and $\mathfrak{a}$.

The second question asks whether  it is possible to create a first countable \PNC\ space with countable extent   in ZFC.

While we could not resolve the first question,  
we made some progress on the second. 
In Section \ref{sc:DRC-notDCC-small}
we prove the following statement (which follows immediately from  Theorem \ref{tm:main-cf}):
\begin{theorem}\label{tm:s}
  If CH holds, then there is a first countable, pseudocompact, but not countably compact space 
  with  $\sss(X)={\omega}$.
\end{theorem}

It is important to emphasize that  ZFC can not  guarantee the existence of such a space , 
as its existence  would imply the existence of an S-space
(see Proposition \ref{pr:S-space}).

 \subsection*{Notions and notations.}

 Given a space $X$ and  a set $A\subs X$
 write 
 \begin{displaymath}
 \accu{A}{X} =\{p\in X: \text{$p$ is an accumulation point of $A$ in $X$}\},
 \end{displaymath} 
and let 
 \begin{displaymath}
 \CD(X)=\{A\in {[X]}^{{\omega}}:\accu{A}{X}=\empt\}=\{A\in {[X]}^{{\omega}}:\text{$A$ is closed discrete} \}.
 \end{displaymath}

 \begin{definition}
 Let $X$ be a topological space and $Y\subs X$. We say that $Y$ is {\em relatively countably compact in $X$},
 and we write $Y\subs^{RC}X$
 iff every infinite subset of $Y$ has an accumulation point in $X$.
 (In \cite{Be81}, Berner referred to this as  "{\em conditionally compact}".)
 
 We write $Y\subs^{DRC}X$ if $Y$ is  both dense and  relatively countably compact in $X$. 
 
 We say that $Y\subs X$ is {\em anti-countably compact (\AC, in short) in $X$} iff ${[Y]}^{{\omega}}\subs \CD(X)$.
 
 \end{definition}

\section{\DRC\ but $\neg$\DCC\ spaces.}\label{sc:DRC-notDCC}

In this section we will construct consistent examples of \Rhered\ ``\DRC\ but $\neg$\DCC'' spaces.
\begin{theorem}\label{tm:drc-notdcc}
  (1)
  If CH holds, then  there is a   crowded 0-dimensional $ T_2$ 
  space $X$ such that  
  \begin{enumerate}[(a)]
  \item   $X$ has a partition $S\cup Y$, where  $S$ is countable and dense, and  $|\overline A|=|X|$ for each  $A\in {[S]}^{{\omega}}$,
\item every $B\in {[Y]}^{{\omega}}$ is closed and discrete in $X$,
\item every countably compact subset of $X$ is  scattered.
  \end{enumerate}
  
  \noindent 
(2) It is consistent that CH holds, $2^{{\omega}_1}$ is as large as you wish, and  
there is a   0-dimensional $ T_2$ 
  space $X$ with $|X|=2^{{\omega}_1}$ such that (a)-(c) above hold for $X$.
\end{theorem}

\begin{proof}[Proof of Theorem \ref{tm:Thm1.2} from Theorem \ref{tm:drc-notdcc}]
  $X$ is \DRC\ as (a) implies that  $S$ is relatively countably compact in $X$.
  Moreover,  since $X$ is crowded, (c) implies that 
  a dense subset of a non-empty regular closed subset $H$  of $X$ can not be countably compact.   
\end{proof}

Before proving Theorem \ref{tm:drc-notdcc}
we need some preparation. 
\begin{definition}\label{df:nice-triple-ch}
 \noindent        
    (1)A triple $\mf X=\<\mc X,\mc B, \mc F\>$ is a {\em nice triple} iff
\begin{enumerate}[(a)]
    \item  $\mc X=\<X,{\tau}\>$ is a  crowded,   0-dimensional space, 
    \item   $X=C\cup \mbb Q $ for some 
 set $C$ of ordinals,
\item  $\mc B=\{B_i:i\in I\}$ is a clopen base of $\mc X$, where 
 $I$ is a set of ordinals  with $|I|=|X|$, 
 \item  the set $\mbb Q$ is dense in $\mc X$, 
 \item $\mc F\subs X\times  {[\mbb Q]}^{{\omega}}$ and $|\mc F|\le |X|$,
\item \label{pr:Fext}
 if $\<a,A\>\in \mc F$, then $a\in \accu{A}{\mc X}$. 
\end{enumerate}
We say that  $\mf X$ is {\em countable} iff 
$X$ is countable.

Observe that we did not assume that the topology ${\tau}$ is $T_2$.

If $\mf X_\ell$ is a nice triple, we will use the notation  $\mc X_\ell$, $X_\ell$,  ${\tau}_{{\ell}}$,    $C_\ell$, $\mc B_\ell$,
$I_\ell$, $B_\ell(i)$ for $i\in I_\ell$, and $\mc F_\ell$.  
    \medskip
    
    \noindent (2)  If 
    $\mf X_0$ and $\mf X_0$  
    are nice triples, then we say 
that   $\mf X_1$ is  an \emph{extension of $\mf X_0$}, and we write  $\mf X_1\ll \mf X_0$, iff
    \begin{enumerate}[(i)]
    \item \label{ext:1x}  $C_0\subs C_1$ and $I_0\subs {I}_1$,  \smallskip
    \item \label{ext:2x} $B_0(i)=B_1(i)\cap X_0$ for each $i\in {I}_0$, \smallskip
    \item \label{ext:3x} if $B_0(i)\subs B_0(i')$ 
    then $B_1(i)\subs B_1(i')$ for each $i,i'\in {I}_0$, \smallskip 
    \item \label{ext:4x} if $B_0(i)\cap B_0(i')=\empt$ 
    then $B_1(i)\cap  B_1(i')=\empt$ for each $i,i'\in {I}_0$, \smallskip
    \item  \label{ext:conv}  $\mc F_0\subs \mc F_1$. 
    \end{enumerate} 
 \end{definition}

    \begin{lemma}\label{lm:limitx}
      Assume that $\<L,\triangleleft\>$ is a directed poset, and 
$\{\mf X_i:i\in L\}$ is a family of countable nice triples such that 
$i\triangleleft j$ implies that $\mf X_j\ll \mf X_i$. 

Then there is a unique nice triple  $\mf X_*$ denoted by 
$\lim_{{\zeta\in L}}\mf X_{\zeta}$, such that 
\begin{enumerate}[({e}1)]
\item $\mf X_*\ll \mf X_{\zeta}$ for each ${\zeta}\in L$,
\item $X_*=\bigcup_{{\zeta}\in L} X_{{\zeta}}$.
\item ${I}_*=\bigcup_{{\zeta}\in L} {I}_{{\zeta}}$.
\item $\mc F_*=\bigcup_{{\zeta}\in L} \mc F_{{\zeta}}$.
\end{enumerate}
If $|L|\le {\omega}$, then $\lim_{{\zeta}\in L}\mc X_{\zeta}$ is countable. 
\end{lemma}

\begin{proof}
Write $C_*=\bigcup_{{\zeta}\in L}{C_{\zeta}}$, $X_*=C_*\cup \mbb Q$,
$I_*=\bigcup_{{\zeta}\in L} {I_{\zeta}}$, 
$\mc F_*=\bigcup_{{\zeta}\in L}\mc F_{{\zeta}}$, 
for  
$i\in I_*$ let 
\begin{displaymath}
B_*(i)=\bigcup\{B_{{\xi}}(i):i\in {I}_{{\xi}}\},
\end{displaymath}
and $\mc B_*=\{B_*(i):i\in I_*\}$. 
Then  $\mc B_*$ is a base of a 0-dimensional topology ${\tau}_*$ on 
$X_*$.
Write $\mc X_*=\<X_*,{\tau}_*\>$.
Then $\mf X_*=\<\mc X_*,\mc B_*,\mc F_*\>$ is a nice triple which meets the requirements, and it is clearly unique.
\end{proof}

\begin{lemma}
  \label{lm:t2}
If $\mf X_0$
is a countable nice triple, 
then there is a countable  extension $\mf X_1$
of $\mf X_0$  such that $X_1=X_0$,   
$\mc X_1$ is $T_2$, and $C_0$ is a closed discrete subspace in $\mc X_1$.
\end{lemma}

\begin{proof}
We can assume that 
$\<x,{\mbb Q}\>\in \mc F_0$ for each $x\in X_0$ because ${\mbb Q}$ is dense in $\mc X_0$.
Consider the family 
\begin{displaymath}
\mc M=\{B_0(i)\cap F: i\in {I}_0, \<{\gamma},F\>\in \mc F_0, {\gamma}\in B_0(i)\}.
\end{displaymath}        
Since  $\mc M\subs {[\mbb Q]}^{{\omega}}$ and $|\mc M|\le {\omega}$,  
 we can choose a family $\mc S=\{S_n:n<{\omega}\}\subs {[{\mbb Q}]}^{{\omega}}$ such that 
\begin{displaymath}
\forall {\varepsilon}\in Fn({\omega},2)\ \forall M\in \mc M\ 
|M\cap S[{\varepsilon}]|={\omega}, 
\end{displaymath}
where $S[\empt]=\mbb Q$, and  
$S[{\varepsilon}]=\bigcap_{{\varepsilon}(n)=1}S_n\cap \bigcap_{{\varepsilon}(n)=0}(\mbb Q\setm S_n)$ for ${\varepsilon}\ne \empt$.

Fix an enumeration  $\big\{\{x_n,y_n\}:n<{\omega}\big\}$ of ${[X_0]}^{2}$, and let 
$$T_n=S_n\cup\{x_n \}\setm \{y_n\}.$$
Consider the family 
\begin{displaymath}
  \mc B'=\{B_0(i)\cap T[{\varepsilon}]:i\in I_0, {\varepsilon}\in Fn({\omega},2)\},
  \end{displaymath}
  where $T[\empt]=X_0$, and  
$T[{\varepsilon}]=\bigcap_{{\varepsilon}(n)=1}T_n\cap \bigcap_{{\varepsilon}(n)=0}(X_0\setm T_n)$ for ${\varepsilon}\ne \empt$.

Then $\mc B'$ is a neighborhood base of a $0$-dimensional  topology ${\tau}_1$ on 
$\mbb Q\cup C_0$.
The topology is $T_2$ because $\{x_n,y_n\}\in {[X_0]}^{2}$ are separated by 
$T[\{\<n,1\>\}]=T_n\ni x_n$ and $T[\{\<n,0\>\}]=X_0\setm T_n\ni y_n$.

The subset $C_0$ is closed discrete, because 
 $x_n\in  T[\{\<n,1\>\}]$ and   $(C_0\setm \{x_n\})\subs T[\{\<n,0\>\}]$.

Moreover, $a\in \accu A{{\tau}_1}$ for each $\<a,A\>\in \mc F_0$.
Indeed, 
 if $\<a,F\>\in \mc F_0$, 
             and $a\in B_0(i)\cap T[{\varepsilon}]$
             then $F\cap B_0(i)$ is infinite  as $a\in B_0(i)$. Since $F\cap B_0(i)\in \mc M$,
             it follows that $F\cap B_0(i) \cap T[{\varepsilon}]$ is also infinite. 
Since  $\<x,\mbb Q\>\in \mc F_0$ for each $x\in X_0$ , ${\tau}_1$ is crowded and $\mbb Q$ is dense in it.

Fix an enumeration $\{B_1(i):i\in I_1\}$ of $\mc B'$ such that 
$B_1(i)=B_0(i)$ for $i\in I_0$.

 Then $\mf X_1=\<\<X_0,{\tau}_1\>,\mc B_1,\mc F_0\>$ meets the requirements.
 \end{proof}

\begin{lemma}\label{lm:extension2}
If $\mf X_0$
is a nice countable triple,  and   
$A\in {[{{\mbb Q}}]}^{{\omega}}$,
then there is a countable  extension 
$\mf X_1$
of $\mf X_0$ such that $X_1=X_0$ and 
$A$ contains an infinite  closed discrete subset $B$  in $\mc X_1$. 
\end{lemma}

    \begin{proof}[Proof of Lemma \ref{lm:extension2}]
      By Lemma \ref{lm:t2}, we can assume that $\mc X_0$ is $T_2$.
      We can also assume that $A$ is not  closed discrete in $\mc X_0$.
Thus, $A$ should contain convergent sequences. 
So we can assume that $A$ converges to some 
${\gamma}$ in $\mc X_0$.

Let $\{B'(\ell):\ell<{\omega}\}$ be an enumeration of $\mc B_0$,
and 
let $\{F_n:n<{\omega}\}$ be an ${\omega}$-abundant enumeration of 
$\{F:\<{\gamma},F\>\in \mc F\}$.

By induction on $n$,  choose $U_n\in \mc B_0$ and $d_n\in A$ such that 
\begin{enumerate}[(i)]
\item $U_n\subs\bigcap\{B'(\ell):\ell<n, {\gamma}\in B'(\ell)\}\setm 
\bigcup\{B'(\ell):\ell<n, {\gamma}\notin B'(\ell)\}$,
\item $U_n\cap F_n\ne \empt$, 
\item ${\gamma}\notin U_n$, $\{d_m:m<n\}\cap U_n=\empt$, 
\item $d_n\in  A\setm \{d_m:m<n\}\setm \bigcup\{U_m:m\le n\}$. 
\end{enumerate}

Let $$V=\{{\gamma}\}\cup \bigcup_{n\in {\omega}}U_{n},$$
and write 
\begin{displaymath}
\mc B_1=\mc B_0\cup\{V\cap B:{\gamma}\in B\in \mc B_0\}.
\end{displaymath}
Then $\mc B_1$ is the neighborhood base of a 0-dimensional  topology 
${\tau}_1$ on $X_0$
such that 
$B=A\setm V$ is an infinite, closed discrete set in ${\tau}_1$.  

By (ii), ${\gamma}\in \accu F{\mc X_1}$ for each $\<{\gamma},F\>\in \mc F_0$.

Fix an enumeration 
 $\{B_1(i):i\in I_1\}$ of $\mc B_1$  such that $B_1(i)=B_0(i)$ for $i\in I_0$.
Then $\mf X_1=\<\<X_0,{\tau}_1\>,\mc B_1,\mc F_0\>$ meets the requirements. 

\end{proof}

\begin{lemma}\label{lm:komb3}

   If $\mf X_0$
   is a nice countable triple,    
   $A\in {[{\mbb Q}]}^{{\omega}}$ is closed discrete in $\mc X_0$,
   and $z\notin C_0$ is an ordinal, 
   then there is a countable extension $\mf X_1$
   of $\mf X_0$ such that $C_1=C_0\cup\{z\}$ and   $\<z,A\>\in \mc F_1$.
   
\end{lemma}

\begin{proof}[Proof of Lemma \ref{lm:komb3}]
 
  We can assume that $\<x,{\mbb Q}\>\in \mc F_0$ for each 
  $a\in X_0$.

        Let $\{B_i:i<{\omega}\}$ be an enumeration of the 
        base $\mc B_0$.
    
        By induction choose a decreasing sequence 
        $\{A_n:n<{\omega}\}$ of infinite subsets of $A$
        such that  
        \begin{displaymath}
        A_n\subset B_n \text{ or } 
        A_n\cap B_n=\empt 
        \end{displaymath} 
for $n<{\omega}$.
Pick pairwise distinct $a_n\in A_n$ for $n\in {\omega}$, then 
choose pairwise disjoint clopen neighborhoods $U_n$ of $a_n$
such that $U_n\subs B_i$ iff $a_n\in B_i$ and
$U_n\cap B_i=\empt$ iff $a_n\notin B_i$ for each $i\le n$.

Then, 
for each $i<{\omega}$, 
\begin{displaymath}
    \forall^\infty n ( U_n\subset B_i)\ \lor\ 
    \forall^\infty n ( U_n\cap  B_i)=\empt.
    \end{displaymath}

Let $C_1=C_0\cup\{z\}$,   
and 
$I_1=I_0\cup\{{\zeta}_n:n<{\omega}\}$, where  ${\zeta}_n\notin I_0$.
For ${\zeta}\in I_0$
let 
\begin{displaymath}
B_1({\zeta})=\left\{\begin{array}{ll}
{B_0({\zeta})}&\text{if $\forall^\infty n ( U_n\cap  B_0({\zeta})=\empt)$},\\\\
{B_0({\zeta})}\cup\{z\}&\text{if $\forall^\infty n ( U_n\subs B_0({\zeta}))$.}
\end{array}\right.
\end{displaymath}

Moreover, for $n<{\omega}$
let 
    \begin{displaymath}
    B_1({\zeta}_n)=\{z\}\cup\bigcup_{m\ge n}U_m. 
    \end{displaymath}
    
  Let ${\tau}_1$ be the topology generated by $\mc B_1=\{B_1(j):j\in I_1\}$ as a base.   
To show that every $B_1(i)$ is closed, assume that $z\notin B_1(j)$. Then 
there is $m\in {\omega}$ such that $U_n\cap B_i(j)=\empt$ for each $n\ge m$. Thus 
$B_1(j)\cap B_1({\zeta}_m)=\empt$.    
  
Finally, put $\mc F_1=\mc F_0\cup\{\<z,A\>\}$. 

Then $\mf X_1=\<\<X_1,{\tau}_1\>,\mc B_1, \mc F_1\>$ satisfies the requirements. 
\end{proof}

\begin{proof}[Proof of Theorem \ref{tm:drc-notdcc}.(1)]
  Let 
  $\<K_0,K_1\>$ be a partition of ${\omega}_1$ into uncountable pieces,  
  and let $\{A_{\xi}:{\xi}\in K_1\}$
  be  an ${{\omega}_1}$-abundant enumeration of the family
  ${[\mbb Q]}^{{\omega}}$.

  We define a $\ll$-decreasing sequence 
  $\<\mf X_{\zeta}:{\zeta}\le {{\omega}_1}\>$
  of nice triples 
  such that 
  \begin{enumerate}[(i)]
  \item   $C_{\zeta}\in {\omega}_1+1$, and  $|X_{\zeta}|=|{\zeta}|+{\omega}$,  
  \item $X_0=\mbb Q$ and ${\tau}_0$ is the usual topology on $\mbb Q$,
  \item if ${\zeta}$ is a limit ordinal,  let 
  $\mf X_{\zeta}=\lim_{{\xi}\in {\zeta}}\mf X_{\xi}$
  (see Lemma \ref{lm:limitx}).
  \item Assume that  ${\zeta}={\xi}+1$, and ${\xi}\in {K_0}$.
  
  Apply Lemma \ref{lm:t2}
  for  $\mf X_{\xi}$ to obtain a countable nice triple 
  $\mf X_{\zeta}$ such that $\mc X_{\zeta}$ is $T_2$ and 
 the countable subset $C_{\xi}$ closed discrete in $\mc X_{{\zeta}}$. 
  \item Assume  ${\zeta}={\xi}+1$, and ${\xi}\in {K_1}$.
  
  First,  apply  Lemma \ref{lm:extension2} 
  for the nice triple $\mf X_{\xi}$ and $A_{\xi}$ to find a countable 
  extension $\mf X'_{\xi}$ of $\mf X_{\xi}$
   such that in $\mf X'_{\xi}$
  the set $A_{\xi}$ contains an infinite  closed discrete set $B_{\xi}$.
  
  Then, applying  Lemma \ref{lm:komb3} 
  for $\mf X'_{\xi}$ and $B$, 
  we can obtain a countable  extension $\mf X_{\zeta}$ of
  $\mf X'_{{\xi}}$ such that $\<a,B\>\in \mc F_{{\zeta}}$.
  We can assume that  $C_{\zeta}=C_{\xi}\cup\{a\}\in {\omega}_1$.
  \end{enumerate}
  
  Finally, $\mc X_{{{\omega}_1}}$ satisfies the requirements.
It is $T_2$ because $\mc X_{\zeta}$ is $T_2$ for cofinally many 
${\zeta}$ and $\mf X_{{\omega}_1}\ll \mf X_{\zeta}$. 
  It is \DRC\ because  $\mbb Q$ is a dense, relatively countably compact subset. 
  
  We also have 
  $\Delta(\mc X_{{\omega}_1})={\omega}_1$.
  Indeed, if $B_i\in \mc B_{{\omega}_1}$, then let $A=B_i\cap \mbb Q$. Then 
  $J=\{{\xi}:A_{\xi}=A\}$ is uncountable,  and 
  for each ${\xi}\in J$ we added a new accumulation point to $A$. But these points are in 
  $B_i$.

  To prove (c)  
assume, for contradiction, that $Z\subs X_{{\omega}_1}$ is a countably compact set that is not scattered.  Then,  there exists  an open set $U$
such that   $T=Z\setm U$ is crowded.  Since $T$ is countably compact,  we must have $|T|\ge {\omega}_1$. Hence,  $T\cap {\omega}_1$ is infinite, which is a contradiction because every infinite countable subset of ${\omega}_1$ is closed discrete.   
 \end{proof}

  \begin{proof}[Proof of Theorem \ref{tm:drc-notdcc}.(2)]
    Assume that $GCH$ holds in the ground model, and let ${\kappa}>{\omega}_1$ 
    be an arbitrarily large regular cardinal. 
    
Consider the poset $\mc P=\<P,\ll\>$, where 
    \begin{displaymath}
    P=\{\mf X_*:
    \text{$\mf X_*$ is a nice triple}, 
    \text{$C_*\cup I_*\in {[{\kappa}]}^{\le {\omega}}$}\}.
    \end{displaymath}

    If $D$ and $E$ are sets of ordinals with $tp(D)=tp(E)$, 
    denote ${\rho}_{D,E}$ the  unique $\in$-preserving bijection between 
         $D$ and $E$. 
    \begin{definition}
    Two conditions $\mf X_0$ and $\mf X_1$ are \emph{twins}
    iff 
    \begin{enumerate}[(1)]
        \item  $tp(C_0)=tp(C_{1})$ and $tp(I_0)=tp(I_1)$, \smallskip
        \item ${\rho}_{C_0,C_1}\restriction C_0\cap C_1=id$, 
        and ${\rho}_{I_0,I_1}\restriction I_0\cap I_1=id$, \smallskip
        \item         for each $i\in I_0$,
        $$B_{1}({\rho}_{I_0,I_1}(i))=(B_0(i)\cap \mbb Q)\cup {\rho}_{C_0,C_1}''
        (B_0(i)\cap {\kappa}),$$
        \item
    $    \mc F_1=\{\<{\rho}(a), A\>:\<a,A\>\in \mc F_0\}.$
        \end{enumerate} 
        \end{definition}
    
    \begin{lemma}\label{lm:amalg2}
    If $\mf X_0$ and $\mf X`_1$ are twins, then they are compatible in $P$.
    \end{lemma}
    
    \begin{proof}
    
    Let $C_2=C_1\cup C_2$, $X_2=\mbb Q\cup C_2, $   $I_2=I_0\cup I_1$, and  for $i\in I_2$ let 
        \begin{displaymath}
        {B_2(i)}=\left\{\begin{array}{ll}
        {B_{0}(i)}\cup B_1({{\rho}_{I_0,I_1}(i)})&\text{if $i\in I_{0}\setm I_1$,}\\\\
        {B_{{1}}(i)}\cup B_0({{\rho}^{-1}_{I_0,I_1}(i)})&\text{if $i\in I_{{1}}\setm I_0$,}\\\\
        {B_{0}(i)\cup B_1(i)}&\text{if $i\in I_0\cap I_1$}.
        \end{array}\right.
        \end{displaymath}
    Then $\{B_2(i):i\in I_2\}$ is a base of a 0-dimensional (but 
    typically not Hausdorff)
    topology ${\tau}_2$  on $X_0\cup X_1$. Moreover, $\mc X_0$ and $\mc X_1$ are subspaces of $\mc X_2$. 
    
    Finally, the triple  
     $\<\<X_2,{\tau}_2\>,\mc B_2, \mc F_0\cup \mc F_1\>\in P$ is a common extension of $\mf X_0$ and $\mf X_1$.
    \end{proof}

    The previous lemma clearly implies the following statement:
    \begin{lemma}\label{lm:w2-cc}
    $P$ satisfies ${\omega}_2$-c.c. 
    \end{lemma}
    
    Since 
    $\mc P$ is ${\sigma}$-closed by Lemma  \ref{lm:limitx}, 
    forcing with $P$ preserves cardinals, and $2^{{\omega}_1}$ in the generic extension will be 
    $((|P|)^{{\omega}_1})^V={\kappa}$.

    Let $\mc G\subs P$ be a generic filter. 
    By Lemma \ref{lm:limitx}, we can consider the nice triple 
    $\mbb X_*=\lim \mc G$.
    By trivial density arguments, we obtain that $X_*=\mbb Q\cup {\kappa}$,
     $I_*={\kappa}$ and $X_*$ is $T_2$ by Lemma \ref{lm:t2}.
    
    So we obtain a 0-dimensional $T_2$ space  
    $\mc X_*$ in $V[\mc G]$. 
We show that $\mc X_*$ satisfies the requirements.

    \begin{lemma}\label{lm:gen2}
    $|\accu A{{\tau}_*}|={\kappa}$ for each $A\in {[\mbb Q]}^{{\omega}}$.
    \end{lemma}
    
    \begin{proof}
      Since $\mc P$ is ${\sigma}$-complete, $A$ is in the ground model.
      Fix ${\delta}<{\kappa}$.
    By applying  Lemma \ref{lm:komb3}, 
    we obtain that  
    \begin{displaymath}
    E_{A,{\delta}}=\{\mf X_0\in P: \<{\gamma},A\>\in \mc F_0 \text{ for some } {\delta}<{\gamma}<{\kappa}\}
    \end{displaymath}
    is dense in $P$.
    Thus, there is $\mf X)\in \mc G\cap E_{A,{\delta}}$.
    Hence, $\accu A{{\tau}_*}\setm {\delta}\ne \empt$.
    Thus, $|\accu A{{\tau}_*}|={\kappa}$. 
    \end{proof}

    \begin{lemma}\label{lm:gen1}
      Every $A\in {[{\kappa}]}^{{\omega}}$ is closed discrete in $\mc X_*$.
      \end{lemma}

      \begin{proof}
        Since $\mc P$ is ${\sigma}$-complete, $A$ is in the ground model.
        Fix $a\in {\kappa}\setm A$.
        By Lemma \ref{lm:t2}, the set 
      \begin{displaymath}
      D_{a,A}=\{\mf X_0\in P: a\notin \accu A{{\tau}_0}\}
      \end{displaymath}
      is dense in $P$.
      Thus, there is $\mf X_0\in \mc G\cap D_{a,A}$.
      Hence, $a\notin \accu A{{\tau}_*}$.
      \end{proof}
      
We can prove (c) in the same manner as in proof of part (1).
Assume, for  contradiction,  that $Z\subs X_*$ is a countably compact set that is not scattered.  Then there exists an open set $U\in {\tau}_*$
      such that   $T=Z\setm U$ is crowded.  Since $T$ is countably compact,  it follows that $|T|\ge {\omega}_1$. However, this implies that $T\cap {\kappa}$ is infinite, which is a contradiction because every infinite countable subset of ${\kappa}$ is closed and discrete.

(a) holds by Lemma \ref{lm:gen2}, and 
    (b) holds by Lemma \ref{lm:gen1}.

    This completes the proof of Theorem \ref{tm:drc-notdcc}(2).
    \end{proof}

\section{Pseudocompact spaces without dense,  relatively countably compact subspaces}
\label{sc:P-notDRC}

Berner's $\Sigma$ is \Rhered ``\SP, but $\neg$\DRC'', but it is not first countable. 
Berner's monster is first countable, but not \Rhered ``$\neg$\DRC''.

In this section, we construct a first countably, \Rhered ``\SP, but $\neg$\DRC'' space
which contains as many countable discrete subsets as possible. 
A pseudocompact, first countable  space cannot be anti-countably compact,
as it must  contain convergent sequences.  
As the next best alternative,  in Corollary \ref{cor:ch}, we construct  
spaces where every 
uncountable subset contains an infinite  closed discrete subset.

\begin{theorem}\label{tm:main}
If  $\mf s=\mf c$, then 
    there is an \SP, crowded, 
      first countable 0-dimensional $T_2$ space $X$ with $\Delta(X)=\mf c$ which is left separated in type $\mf c$.
\end{theorem}

Assuming CH, we can get a bit more.

\begin{corollary}\label{cor:ch}
If CH holds, then there is a selectively pseudocompact, crowded, 
first countable 0-dimensional $T_2$ space $X$ with $\Delta(X)={\omega}_1$
such that  relatively countable compact subset is countable,
and every countable set is nowhere dense.
\end{corollary}

To obtain Theorem \ref{tm:Tm1.3} and Corollary \ref{cor:ch}
from Theorem \ref{tm:main},
we need to formulate some results which excludes the existence of 
certain relatively countably compact subspaces in certain left separated spaces.

\begin{theorem}\label{tm:geju}
(1)    A left separated, crowded   regular space $Y$ is not \DRC. 

(2)    A first countable 0-dimensional $T_2$ space   which is  
    left separated in type ${\omega}_1$ is not \DRC.
   \end{theorem}
  
 Observe that in (2) we do not assume that the space is crowded.

\begin{proof}[Proof of Theorem  \ref{tm:Tm1.3} from Theorem \ref{tm:main}
  and Theorem \ref{tm:geju}.(1)]
  Consider the space $X$ we obtain from Theorem \ref{tm:main}. 
If $H$ is a regular closed subset of $X$, then $H$ is crowed,
so it is not \DRC\ by Theorem \ref{tm:geju}.(1).
\end{proof}

\begin{proof}[Proof of Corollary \ref{cor:ch}
  from Theorem \ref{tm:main}
  and Theorem \ref{tm:geju}.(2)]
  Consider the space $X$ we obtain from Theorem \ref{tm:main}.
  Then $X$ is left-separated in type $\mf c={\omega}_1$.
Let $Y$ be an uncountable subset of $X$. Then $Y$ 
is also left-separated in type ${\omega}_1$, and so it
is not 
\DRC\  by Theorem \ref{tm:geju}.(2).
\end{proof}

    \begin{proof}[Proof of Theorem \ref{tm:geju}(1)]
Let $\{y_{\alpha}:{\alpha}<{\kappa}\}$ be a left-separating enumeration of $Y$
and let $D\subs Y$ be dense.  
        
By recursion on $n\in {\omega}$    pick $y_{{\alpha}_n}\in D$ and $U_n, V_n\in {\tau}_Y$ as follows.
    
    Let $y_{{\alpha}_0}\in D$ be arbitrary. 
    
    If $y_{{\alpha}_n}$ is given, let $U_n$ be a left-separating neighborhood of 
    $y_{{\alpha}_n}$. Since $Y$ is  regular,  
    we can choose  $V_n\in {\tau}_Y^+$ such that $\overline{V_n}\subs U_n $.
    
Since $Y$ is crowded, we can     pick $y_{{\alpha}_{n+1}}\in D\cap (V_n\setm \{y_{{\alpha}_n}\})$.

    We claim that  $\{y_{{\alpha}_n}:n<{\omega}\}\in {[D]}^{{\omega}}$ is closed discrete in $Y$.
Indeed, ${\alpha}_n<{\alpha}_{n+1}$ by the construction. Let ${\alpha}=\sup\{{\alpha}_n:n<{\omega}\}$.
Then 
\begin{multline*}
\{y_{{\alpha}_n}:n<{\omega}\}'\subset \bigcap_{n<{\omega}}\overline{\{y_{{\alpha}_m}:m\ge n\}}\subset 
\overline{\{y_{\zeta}:{\zeta}<{\alpha}\}}\cap 
\bigcap_{n<{\omega}}\overline{V_n}\subs\\ 
\{y_{\zeta}:{\zeta}<{\alpha}\}\cap
\bigcap_{n<{\omega}}\{y_{\zeta}:{\alpha}_n\le {\zeta}\}=
\{y_{\zeta}:{\zeta}<{\alpha}\}\cap \{y_{\zeta}:{\alpha}\le {\zeta}\}=\empt. 
\end{multline*}
    \end{proof}

\begin{proof}[Proof of Theorem \ref{tm:geju}(2)] We can assume that $Y={\omega}_1$.  Let $D\subs Y$ be dense. 
  Let $\{B({\alpha},i):{\alpha}<{\omega}_1,i<{\omega}\}$
  be a clopen base of $Y$ such that $B({\alpha},i)\supset B({\alpha},i+1)$ and $B({\alpha},0)\cap {\alpha}=\empt$.

 By induction on $n$ pick ${\alpha}_n\in D$, $ {\beta}_n\in Y$  and $k_n,i_n\in {\omega}$
 such that 
 \begin{enumerate}[(1)]
 \item ${\alpha}_{n-1}<{\alpha}_n$,
 \item $\{{\beta}_n:n<{\omega}\}=\bigcup\{{\alpha}_n:n<{\omega}\}$,
 \item $Y_n=Y\setm\bigcup\{B({\beta}_m,i_m):m\le n\}$ is uncountable,
 \item ${\alpha}_n\in Y_n\cap D$.
 \end{enumerate}
 
 Assume that we have ${\alpha}_m,{\beta}_m, k_m,i_m$ for $m<n$.
 
 Using a bookkeeping function choose ${\beta}_n$ such that (2) will hold.
 Since $Y_{n-1}$ is uncountable, we can choose $i_n$ such that $Y_n=Y_{n-1}\setm B({\beta}_n,i_n)$  
 is still uncountable.

 Since $D$ is dense, $Y_n$ is uncountable clopen, and $Y$ is left-separated,   we can pick ${\alpha}_n\in D\cap (Y_n\setm \max({{\alpha}}_{n-1}+1,{\beta}_{n}+1)$.

 Let ${\delta}=\bigcup\{{\alpha}_n:n<{\omega}\}=\{{\beta}_n:n<{\omega}\}$.
   Then $\overline {\{{\alpha}_n:n<{\omega}\}}\subs {\delta}$ because 
  $Y$ is left-separated, and 
 ${\beta}_m$ is not an accumulation point of $\{{\alpha}_n:n<{\omega}\}$
 because 
 \begin{displaymath}
 B({\beta}_m,i_m)\cap \{{\alpha}_n:n<{\omega}\}\subs \{{\alpha}_k:k\le m\}.
 \end{displaymath}
 Thus, $\{{\alpha}_n:n<{\omega}\}\subs D$ is closed discrete  
 in $Y$. So  $D$ is not relatively countably compact. 
 \end{proof}

 Before proving Theorem \ref{tm:main}, we need to prove some lemmas.

\begin{definition}\label{df:triple-ma}
    (1) A triple $\mf X=\<\mc X,\mc B,\prec_X\>$ is a {\em good triple} 
    iff
    \begin{enumerate}[(t1)]
    \item 
    $\mc X=\<X,{\tau}\>$ is a left-separated, crowded,  first countable,  0-dimensional $T_2$-space,
    \item $\mc B=\<B(x,i):x\in X, i\in {\omega}\>$ is a family of clopen sets, 
    \item 
    $\{B(x,i):i\in {\omega}\}$ is a neighborhood base at $x$ in  $X$ consisting of clopen subsets
    such that $B(x,i)\supset B(x,i+1)$ for each $i<{\omega}$. 

  \item $\prec_X$ is a left separating well-ordering of $X$,
  \end{enumerate}

   \medskip 
If $\mf X_\ell$ is good triple, write 
$\mf X_\ell=\<\mc X_\ell, \mc B_\ell, \prec_\ell\>$,
$\mc X_\ell=\<X_\ell, {\tau}_\ell\>$, moreover let  
$\mc B_\ell=\<B_\ell(x,i):x\in X_\ell, i<{\omega}\>$.

    \medskip
    
    \noindent (2)  Given good triples $\mf X_\ell=\<\mc X_\ell,\mc B_\ell, \prec_\ell\>$
     for $\ell \in 2$,  we say 
that   $\mf X_1$ is  an \emph{extension of $\mf X_0$}, and we write  $\mf X_1\ll \mf X_0$, iff
    \begin{enumerate}[({e}1)]
    \item \label{ext:1}  $X_0\subs X_1$,  \smallskip
    \item \label{ext:2} $B_0(x,i)=B_1(x,i)\cap X_0$ for each $x\in X_0$ and $i\in {\omega}$, \smallskip
    \item \label{ext:3} if $B_0(x,i)\subs B_0(x',i')$ 
    and $x'\notin B_0(x,i)$
    then $B_1(x,i)\subs B_1(x',i')$ for each $x,x'\in X_0$ and $i,i'<{\omega}$, \smallskip 
    \item \label{ext:4} if $B_0(x,i)\cap B_0(x',i')=\empt$ 
    then $B_1(x,i)\cap B_1(x',i')=\empt$ for each $x,x'\in X_0$ and $i,i'<{\omega}$, \smallskip 
    \item \label{ext:5} $\prec_0\subs \prec_1$ and $X_0$ is an initial segment in $\<X_1,\prec_1\>$.
    \end{enumerate} 
    \end{definition}

\begin{klemma}\label{lm:double-ma}
Assume that 
\begin{enumerate}[(a)]
  \item $\mc X_0$ is a  good triple,
  \item  $|X_0|<\mf s$,
 \item the family $\{B_0({\zeta},j({\zeta})):{\zeta}\in K\}$ is locally finite in $X$
for some 
$K\in {[X]}^{{\omega}}$ and $j:K\to {\omega}$, 
\end{enumerate}
Then  there is a good triple $\mf X_1$ 
such that 
\begin{enumerate}[(1)]
\item $\mf X_1\ll \mf X_0 $,
\item the family $\{B_1({\zeta},j({\zeta})):{\zeta}\in K\}$ is not locally finite in $\mc X_1$. 
\item $|X_1|=|X_0|$.
\end{enumerate}
\end{klemma}

\begin{proof}[Proof of the Key Lemma \ref{lm:double-ma}]

For ${\zeta}\in K$ pick ${\eta}_{\zeta}\in B_0({\zeta},j({\zeta}))\setm \{{\zeta}\}$.  
Let $K_*=\{{\eta}_{\zeta}:{\zeta}\in K\}$.
Since $|X_0|<\mf s$, the family $\{B_0(x,i)\cap K_*:x\in X_0,i\in{\omega}\}$
can not be a splitting family on ${[K_*]}^{{\omega}}$. So,   
there is a set $L\in {[K]}^{{\omega}}$
such that writing $L_*=\{{\eta}_{\zeta}:{\zeta}\in K_*\}$ 
 for each $\<x,i\>\in X\times {\omega}$ we have 
\begin{displaymath}
    L_*\subs ^*  B_0(x,i)\   \lor\ 
    L_*\subs^* X_0\setm  B_0(x,i) .  
    \end{displaymath}

The underlying set of the extension  $\mf X_1$ will be 
\begin{displaymath}
X_1=X_0\cup\{p\}\cup (X_0\times \mbb Q),
\end{displaymath}
where $p$ is a new element.

For $q\in \mbb Q$
let $\{I(q,i):i\in {\omega}\}$ be a clopen neighborhood base  of $q$ in $\mbb Q$.
Fix an enumeration $\{{\zeta}_n:n<{\omega}\}$  of $L$.

Define  $B_1(y,i)$ for $y\in X_1$ and $i<{\omega}$ as follows.

\smallskip

\noindent
{\bf Case 1.} $y=\<x,q\>\in X_0\times \mbb Q$. 

Let 
\begin{displaymath}
B_1(y,i)=\{x\}\times I(q,i).
\end{displaymath}

\smallskip

\noindent
{\bf Case 2}. $y=p$.

Let 
\begin{displaymath}
B_1(p,i)=\{p\}\cup\bigcup_{n\ge i}\big (\{{\eta}_{\zeta_n}\}\times \mbb Q\big ).
\end{displaymath}

\smallskip

\noindent
{\bf Case 3}. $y\in X$.

Let 
\begin{displaymath}
B'(y,i)=B_0(y,i)\cup \Big(B_0(y,i)\setm \{y\}\Big)\times \mbb Q,
\end{displaymath}
and 
\begin{displaymath}
B_1(y,i)=
\left\{\begin{array}{ll}
{B'(y,i)}&\text{if $L_*\subs^* X_0\setm B_0(y,i)$,}\\\\
{B'(y,i)}\cup \{p\}&
\text{if $L_*\subs^* B_0(y,i)$.}
\end{array}\right.
\end{displaymath}

 Finally, let $\prec_Q$ be a well-ordering of $\mbb Q$
in type ${\omega}$, and define $\prec_1$ as follows.
\begin{enumerate}[(a)]
\item $\prec_0\subs \prec_1$,
\item $\forall x\in X_0$\  $x\prec_1p$,
\item $\forall y\in X_0\times \mbb Q$  $p\preceq_1 y$,
\item $\prec_1\restriction X_0\times \mbb Q$ is the lexicographical product of 
$\prec_0$ and $\prec_Q$.
\end{enumerate} 

In that way we defined $\mf X_1$. 
We should check first that $\mf X_1$ is a
good triple extending $\mf X_0$.

\begin{claim}
  (e\ref{ext:1}),    (e\ref{ext:2}), (e\ref{ext:4}) and 
     (e\ref{ext:5})   hold for $\mf X_0$ and $\mf X_1$.
\end{claim}

Trivial from definition.

\begin{claim}\label{cl:bsubB'}
    If 
    $B_0(x,i)\subs B_0(x',i')$ and  $p\in B_1(x,i)$, then  $p\in B_1(x',i')$.  
\end{claim}

Indeed, 
 if $L_*\subs^* B_0(x,i)$,
then  $L_*\subs^* B_0(x',i')$.

\begin{claim}\label{cl:e3}
(e\ref{ext:3}) holds for $\mf X_0$ and $\mf X_1$.
\end{claim}

\begin{proof}

Indeed, 
if $B_0(x,i)\subs B_0(x',i')$ 
and $x'\notin B_0(x,i)$
then 
\begin{multline*}
B_1(x,i)\cap (X_0\times \mbb Q)= (B_0(x,i)\setm \{x\})\times \mbb Q
\subs \\(B_0(x',i')\setm \{x'\})\times \mbb Q=
B_1(x',i')\cap (X_0\times \mbb Q),
\end{multline*}
and  $p\in B_1(x,i)$ implies $p\in B_1(x',i)$ by Claim \ref{cl:bsubB'}. 
\end{proof}

\begin{claim}
$\big\{\{B_1(y,n):n\in {\omega}\}:y\in X_1\big\}$ is a neighborhood system of a topology ${\tau}_1$ on $X_1$. 
\end{claim}

\begin{proof}[Proof of the Claim]
    By \cite[Proposition 1.2.3]{en}, we should check that 
    \begin{enumerate}[(BP1)]
    \item     $y\in B_1(y,n)$ for each  $y\in Y$  and $ n<{\omega}$ ,
    \smallskip
    \item if $z\in B_1(y,n)$ then $B_1(z,m)\subs B_1(y,n)$ for some 
    $m<{\omega}$,    \smallskip

    \item for each $x\in X_1$ and for each   $n,m<{\omega}$
    there is $k<{\omega}$ such that $B_1(x,k)\subs B_1(x,n)\cap B_1(x,m) $.
    \end{enumerate}
Conditions (BP1) and (BP3) are trivial. 

To check (BP2), assume that $z\in B_1(y,n)$, $z\ne y $. 
If $y=\<x,q\>\in X_0\times \mbb Q$,
then $z=\<x,r\>$ for some $r\in I(q,n)$. 
Thus, there is $m$ with $I(r,m)\subs I(q,n)$ and so 
 $B_1(z,m)\subs B_1(y,n)$. 

If $y=p$, then $z=\<x,r\>$, where    $x ={\eta}_{{\zeta}_k}$ 
for some $k\ge n$, and $r\in \mbb Q$, 
and so $B_1(z,m)\subs \{{\eta}_{{\zeta}_k}\}\times \mbb Q\subs B_1(p,n)$ for each $m\in {\omega}$.

Finally, consider the case  $y\in X_0$. If $z\in X_0$, then pick $m$ such that 
$B_0(z,m)\subs B_0(y,n)\setm \{y\}$. Then $B_1(z,m)\subs B_1(y,n)$ by 
Claim \ref{cl:e3}.

If $z=\<x',q\>\in X_0\times \mbb Q$, then $x'\ne x$ by the definition of 
$B_1(x,n)$.
Thus,  $B_1(z,m)\subs \{x'\}\times \mbb Q\subs B_1(x,n)$ for each $m\in {\omega}$.

 Now, assume that $z=p$. Then 
there is $m\in {\omega}$ such that 
${\eta}_{{\zeta}_k}\in  B_0(x,n)$
for each $k\ge m$. 
Hence $B_1(p,m)\subs B_1(x,n)$.
\end{proof}

\begin{claim}
${\tau}_1$ is $T_2$.
\end{claim}

\begin{proof}
Fix  $\{y,z\}\in {[X_1]}^{2}$.

Assume first that  $y\in X_0$ and $z=p$. 
Since the family $\{B_0({\zeta},j({\zeta})):{\zeta}\in K\}$ is locally finite in $\mc X_0$,
there are $i,m\in {\omega}$ such that 
$B_0({\zeta}_n,g({\zeta}_n))\cap B_0(y,i)=\empt$ for each $n\ge m$. In particular, ${\eta}_{\zeta_n}\notin B_0(y,i)$ for $n\ge m$,  and so
$B_1(p,m)\cap B_1(y,i)=\empt$.

\smallskip
If $y\in X_0$ and $z=\<x,q\>\in (X_0\setm \{y\})\times \mbb Q$ then 
pick $i$ such that $x\notin B_0(y,i)$. Then 
$B_1(y,i)\cap B_1(z,j)\subs B_1(y,i)\cap (\{x\}\times \mbb Q)=\empt$
for each $j\in {\omega}$.

If $y\in X_0$ and $z=\<y,q\>\in \{y\}\times \mbb Q$ then 
$B_1(y,i)\cap B_1(z,k)\subs B_1(y,i)\cap (\{y\}\times \mbb Q)=\empt$
for each $i, k\in {\omega}$.

The remaining  cases are trivial. 
\end{proof}

\begin{claim}
Every $B_1(y,i)$ is closed, so ${\tau}_1$ is zero-dimensional.  
\end{claim}

\begin{proof}
Fix $z\in X_1$ with $z\notin B_1(y,i)$.

If $\{y,z\}\in {[X_0]}^{2}$, then $z\notin B_0(y,i)$, so we can pick $k$
such that  $B_0(z,k)\cap B_0(y,i)=\empt.$  Then 
$B_1(z,k)\cap B_1(y,i)=\empt$ by (e\ref{ext:4}).

Since 
$X_0\times \mbb Q$ is an open subspace in $\mc X_1$
and 
the subspace topology on $X_0\times \mbb Q$ is the product topology of 
the discrete topology on $X_0$ and the topology of $\mbb Q$, it follows that 
if $\{y,z\}\in {[X_0\times \mbb Q]}^{2}$ then there is $k$ such that 
$B_1(z,k)\cap B_1(y,i)=\empt.$

Consider next the case when   $y\in X_0$ and $z=\<x,q\>\in X_0\times \mbb Q$. 
Then   $z\notin B_0(y,i)$ implies  $B_1(y,i)\cap (\{x\}\times \mbb Q)=\empt$ and so    
$B_1(y,i)\cap B_1(z,k)=\empt$ for each $k\in {\omega}$. 

Assume next that $y\in X_0$ and $z=p$. 
Then $z\notin B_1(y,i) $ implies that   there is  $m\in {\omega}$ such that
${\eta}_n\notin  B_0(y,i)$
for each $n\ge m$.
Thus, $B_1(y,i)\cap B_1(p,m)=\empt$.

 Finally, assume that $y=p$. 

Consider first that case  $z\in X_0$. 
Since  the family $\{B_0({\zeta},j({\zeta})):{\zeta}\in K\}$ is locally 
finite, there is 
$m$ such that 
$K\cap B_0(z,m)$ contains at most one element, namely $z$.
Then,  $B_1(p,i)\cap B_1(z,m)=\empt.$

 Now, assume that $z=\<x,q\>\in X\times \mbb Q$. Then 
$x\notin \{{\eta}_{{\zeta}_\ell}: i\le \ell<{\omega}\}$, so 
$B_1(p,i)\cap B_1(z,j)=\empt$ for each $j\in {\omega}$.
\end{proof}

\begin{claim}
$\prec_1$ is a left-separating well ordering of $X_1$.
\end{claim}
\begin{proof}
    Trivial. 
\end{proof}

Putting together these observations we obtain that 
$$\mf X_1=\<\mc X_1, \{B_1(y,i):y\in X_1, i\in {\omega}\},\prec_1\>$$ is a good triple and 
$\mf X_1\ll \mf X_0$.  Moreover, $p$ is an accumulation point of the family $\{B_1({\zeta},j({\zeta})):{\zeta}\in K\}$.
\end{proof}

\begin{lemma}\label{lm:limit}
Assume that $\<I,\triangleleft\>$ is a directed poset, and 
$\{\mf X_i:i\in I\}$ is a family of good triples such that 
$i\triangleleft j$ implies that $\mf X_j\ll \mf X_i$. 
Then there is a  good triple $\<\mc X_*,\mc B_*, \prec_*\>$   denoted by 
$\lim_{i\in I}\mf X_i$, such that 
\begin{enumerate}[(a)]
\item $\lim_{i\in I}\mf X_i\ll \mf X_i$ for each $i\in I$,
\item $X_*=\bigcup_{i\in I} X_{i}$.
\end{enumerate}
\end{lemma}

\begin{proof}
Write $X_*=\bigcup_{i\in I}X_i$,  and for 
$x\in X_*$ and for $n\in {\omega}$ let 
\begin{displaymath}
B_*(x,n)=\bigcup\{B_j(x,n):x\in X_j\},
\end{displaymath}
and put
\begin{displaymath}
\prec_*=\bigcup_{i\in I}\prec_{i}.
\end{displaymath}
Then, $\mc B_*$ is a base of a 0-dimensional $T_2$ topology ${\tau}_*$ on $X_*$.
Thus,  
writing $\mc X_*=\<X_*,{\tau}_*\>$ the triple
$\mc Z_*=\<\mc X_*, \{B_*(x,n):x\in X_*, n\in {\omega}\},\prec_*\>$ 
satisfies the requirements.
\end{proof}

\begin{proof}[Proof of Theorem \ref{tm:main}]
    
         Let $\{\<K_{\alpha},j_{\alpha}\>:{\alpha}<{\mf c}\}$
         be a ${\mf c}$-abundant enumeration of the family
         \begin{displaymath}
         \{\<K,j\>:K\in {[{\mf c}]}^{{\omega}}, j:K\to {\omega}\}.
         \end{displaymath}
    
         We define a $\ll$-decreasing sequence 
         $\<\mf X_{\zeta}:{\zeta}\le {\mf c}\>$
         of good triples
     such that 
    \begin{enumerate}[(i)]
      \item   $X_{\zeta}$ is an ordinal,  $|X_{\zeta}|=|{\zeta}|+{\omega}$, 
 and      $\prec_{\zeta}$ is the natural ordering of ordinals, 
      \item $X_0$ is a crowded 0-dimensional, first countable  $T_2$
    topology on ${\omega}$,
    \item if ${\zeta}$ is a limit ordinal,  let 
    $\mf X_{\zeta}=\lim_{{\xi}\in {\zeta}}\mf X_{\xi}$
    (see Lemma \ref{lm:limit}).
    \item If ${\zeta}={\xi}+1$, do the following. 
    \begin{enumerate}[(a)]
    \item 
    Consider $K_{\xi}$ and $j_{\xi}$. 
\item    
    If $K_{\xi}\notin{[X_{\xi}]}^{{\omega}}$
    or $\{B_{\xi}(k,j_{\xi}(k)):k\in K_{\xi}\}$ is not a locally finite family of open sets in $\mc X_{\xi}$, 
    then we do nothing, i.e. let $\mf X_{\zeta}=\mf X_{\xi}$.

    \item 
    If 
    $K_{\xi}\in{[X_{\xi}]}^{{\omega}}$
    and  $\{B_{\xi}(k,j_{\xi}(k)):k\in K_{\xi}\}$ is a locally finite family of open sets in $\mc X_{\xi}$,
    apply Lemma \ref{lm:double-ma} for $\mf X_{\xi}$ and 
    $\{B_{\xi}(k,j_{\xi}(k)):k\in K_{\xi}\}$
    to obtain $\mf X_{\zeta}$. Hence,  $\{B_{{\xi}+1}(k,j_{\xi}(k)):k\in K_{\xi}\}$ is not locally finite. 
    We can assume that $X_{\zeta}\in \mf c$ is an ordinal, and $\prec_{\zeta}$ is the natural ordering on 
    that ordinal.  
    \end{enumerate}
\end{enumerate}
    
    Finally, $X_{{\mf c}}$ satisfies the requirements.
To show that $\mc X_{{\mf c}}$ is \SP, let $\{B_\mf c(k,j(k)):k\in K\}$
be a family of basic open sets. 
There is ${\xi}<\mf c$ such that $K_{\xi}=K\in {[X_{\xi}]}^{{\omega}}$ and $j_{\xi}=j$.
Then, by  the construction,  
$\{B_{\xi+1}(k,j(k)):k\in K\}$ is not locally finite, it has an accumulation point $p$.
Since $\chi(p,\mc X_{\xi+1})={\omega}$, we can pick
$x_k\in B_{\xi}(k,j(k))$ for $k\in K$
such that $p \in \accu{\{x_k:k\in K\}}{ {\tau}_{\xi+1}}$.  
Since   $\mf X_{\mf c }\ll \mf X_{\xi+1}$, we have ${\tau}_{\xi+1}=\{U\cap X_{\xi+1}:U\in {\tau}_{\mf c}\}$. 
Hence,  $p \in \accu{\{x_k:k\in K\}}{ {\tau}_{\mf c}}$.
\end{proof}

 \section{A pseudocompact, but not countably compact space with countable spread}
\label{sc:DRC-notDCC-small}

 First,  we make  the following observation:
 the proposition below implies that 
 ZFC alone is insufficient to construct a space as required in Theorem \ref{tm:s}.

     \begin{proposition}\label{pr:S-space}
     If there is a pseudocompact, but not countably compact, regular space $X$  
     with $s(X)={\omega}$, then there is an $S$-space.
     \end{proposition}
 
 \begin{proof}
 A Lindelöf pseudocompact space is compact, so $X$ can not be 
 Lindelöf, and so  it contains a right-separated subspace 
 $Y\in {[X]}^{{\omega}_1}$. Since $s(Y)\le s(X)={\omega}$,
 it follows that $z(Y)={\omega}$ as well. Thus,  $Y$ is an S-space.
 \end{proof}

\begin{theorem}\label{tm:main-cf}
  If $CH$ holds, then there is  a \DRC, but $\neg$\DCC, locally countable, 
  locally compact, first countable, 0-dimensional $T_2$ space $X$ with cardinality ${\omega}_1$ 
  and     $\sss(X)={\omega}$.   
      \end{theorem}

We do not know if we can find an R-hereditary example for the problem we addressed in the previous theorem.

 Before proving Theorem \ref{tm:main-cf}    we need some preparation.
   The first statement is well-known:
   \begin{lemma}\label{lm:known}
   If $Y$ is a countable, regular space, $D\subs Y$ is closed discrete,
   then there is a neighborhood assignment $W:D\to {\tau}_Y$
   such that the family $\{W(d):d\in D\}$ is closed discrete. 
   \end{lemma}

 The Euclidean topology on $\mbb R$ is denoted by ${\varepsilon}$.
 The next lemma is the key of our proof. 
 
  \begin{lemma}\label{lm:induction-lemma}
  Assume that 
  \begin{enumerate}[(a)]
     \item $Y=\<Y,{\tau}_Y\>$  is a countable, locally compact $T_2$ space,  
     \item $Y\cap \mbb R$ is closed in $Y$ and    $p\in \mbb R\setm Y$,
 \item the topology ${\tau}_Y\restriction Y\cap \mbb R$ 
 refines the Euclidean topology on 
 $Y\cap \mbb R$,
 \item $E\in {[Y\setm \mbb R]}^{{\omega}}$ is closed discrete in $Y$,
 \item $\mc D\subs {[Y\cap \mbb R]}^{{\omega}}$, $\mc D$ is countable.   
 \end{enumerate}
  Then there is a  
  space $Z=\<Z,{\tau}_Z\>$ such that 
   \begin{enumerate}[(a')]
  \item $Z$ is locally compact $T_2$,
     \item $Z=Y\cup\{p\}$ and ${\tau}_Y={\tau}_Z\cap \mc P(Y)$,
     \item the topology ${\tau}_Z\restriction Z\cap \mbb R$ refines 
 the Euclidean topology on  
     $Z\cap \mbb R$,
     \item $p\in \overline{E}^Z$,
   \item for each $D\in \mc D$ if  $p\in \overline{D}^{{\varepsilon}}$, then 
   $p\in \overline{D}^{Z}$.
   \end{enumerate}
  \end{lemma}

 \begin{proof}[Proof of Lemma \ref{lm:induction-lemma}]
     Write $S=Y\cap \mbb R$ and  $A=Y\setm \mbb R$. 
     Fix an enumeration $E=\{e_n:n\in {\omega}\}$.
 
  \noindent {\bf Case 1.} {\em $p\notin \overline{S}^{{\varepsilon}}$.}   
 
  By Lemma \ref{lm:known}, there is 
 a neighborhood assignment $W:E\to {\tau}_Y$
 such that the family $\{W(e):e\in E\}$ is closed discrete.
  For each $k\in {\omega}$ write
  \begin{displaymath}
  V_k=\{p\}\cup\{W(e_n):n\ge k\},   
 \end{displaymath}
 and define the topology of ${\tau}_Z$ as follows:
 \begin{enumerate}[(a)]
 \item $\<Y,{\tau}_Y\>$ is an open subspace of $\<Z,{\tau}_Z\>$,
 \item $\{V_k:k\in {\omega}\}$ is a neighborhood base of 
 $p$ in $\<Z,{\tau}_Z\>$.
 \end{enumerate} 
 Then $\<Z,{\tau}_Z\>$ clearly satisfies the requirements.  
 
 \medskip
 \noindent {\bf Case 2.} {$p\in \overline{S}^{{\varepsilon}}$.}
 
  Choose a sequence $P=\{p_n:n<{\omega}\}\subs S$
 such that $$\lim_{\varepsilon}{\{p_n:n<{\omega}\}}= p,$$ and 
 for each $D\in \mc D$, if $p\in \overline{D}^{{\varepsilon}}$, then 
 $D\cap \{p_n:n<{\omega}\}$ is infinite. 
 
 Let $U_n$ be a compact open neighborhood of $p_n$ in ${\tau}_Y$ for $n\in {\omega}$
 such that the family $\{U_n\cap \mbb R:n<{\omega}\}$ 
 converges to $p$ in the 
 Euclidean topology.
 
 Since $S$ is closed in $Y$, we have that $P\cup E$ is closed discrete in 
 $Y$. Thus, by Lemma \ref{lm:known}, there is 
 a neighborhood assignment $W:P\cup E\to {\tau}_Y$
 such that the family $\{W(x):x\in P\cup E\}$ is closed discrete. 
 We can assume that $W(p_n)\subs U_n$.

 For $k\in {\omega}$
 write 
 \begin{displaymath}
 V_k=\{p\}\cup\bigcup_{n\ge k}(W(p_n)\cup W(e_n)).
 \end{displaymath}
 Define the topology of ${\tau}_Z$ as follows: 
 \begin{enumerate}[(a)]
 \item $\<Y,{\tau}_Y\>$ is an open subspace of $\<Z,{\tau}_Z\>$,
 \item $\{V_k:k\in {\omega}\}$ is a neighborhood base of 
 $p$ in $\<Z,{\tau}_Z\>$.
 \end{enumerate}
 This construction clearly works. 
 \end{proof}
 
 \begin{proof}[Proof of Theorem \ref{tm:main-cf}]
 Let
     $\{D_{\zeta}:{\zeta}<{\omega}_1\}={[\mbb R]}^{{\omega}}$,
     $\{p_{\xi}:{\xi}<{\omega}_1\}=\mbb R$, and 
     $\{E_{\xi}:{\xi}<{\omega}_1\}={[{\omega}\times {\omega}]}^{{\omega}}$.
 
 We will define a sequence $\<\<X_{\alpha},{\tau}_{\alpha}\>:{\alpha}\le {\omega}_1\>$
 of countable, locally compact  $T_2$ spaces such that 
 \begin{enumerate}[(a)]
 \item $X_{\alpha}=({\omega}\times ({\omega}+1))\cup\{p_{\zeta}:{\zeta}<{\alpha}\}$,
 \item ${\tau}_{\beta}\cap \mc P(X_{\alpha})={\tau}_{\alpha}$ for 
 ${\alpha}<{\beta}$,
 \item ${\omega}\times {\omega}$ is dense in ${\tau}_{\alpha}$,
 \item $E_{\alpha}$ has an accumulation point in ${\tau}_{{\alpha}+1}$,
 \item if ${\zeta}<{\alpha}$ and   $D_{\zeta}\subs X_{\alpha}$ and 
 $p_{\alpha}\in \overline{D_{\zeta}}^{{\varepsilon}}$, then 
 $p_{\alpha}\in \overline{D_{\zeta}}^{{\tau}_{{\alpha}+1}}$.
 \end{enumerate}
 
 We have  $X_0=({\omega}\times ({\omega}+1))$, and 
     let topology  ${\tau}_0$ on $({\omega}\times ({\omega}+1))$ be  the product topology.
 
     In limit step, take the direct limit.

 To get $X_{{\alpha}+1}$ from $X_{\alpha}$ apply 
 Lemma \ref{lm:induction-lemma} for $Y=X_{\alpha}$, 
 $\mc D=\{D_{\xi}:{\xi}<{\alpha}\}$,  $E=E_{{\beta}(\alpha)}$  and  $p=p_{\alpha}$,
 where $${\beta}({\alpha})=\min\{{\beta}:\text{$E_{\beta}$ is closed discrete in ${\tau}_{\alpha}$}\}.$$
 
 The space $\mc X=\<X_{{\omega}_1},{\tau}_{{\omega}_1}\>$ is clearly locally countable, locally compact, 0-dimensional $T_2$
 with cardinality ${\omega}_1$. The subspace ${\omega}\times {\omega}$ is 
 dense and relatively countably compact because every $E_{\alpha}$ has accumulation point, 
 so $X$ is \DRC. 
 
 If $D\subs X_{{\omega}_1}$ is dense, then $D$ should contain 
 the isolated point: ${\omega}\times {\omega}\subs D$.
Since $\{n\}\times {\omega}$ converges to $\<n,{\omega}\>$,
if $D$ is countable compact, then $E=\{\<n,{\omega}\>:n<{\omega}\}\subs D$.
But $E$ is closed discrete in $\mc X$, so $D$ can not be countably compact.
Thus, $\mc X$ is $\neg$\DCC.

Finally,  if $D\in {[\mbb R]}^{{\omega}_1}$, then $D$ has a countable ${\varepsilon}$-dense subset 
 $D_{\zeta}$. Pick ${\alpha}>{\zeta}$ such that $p_{\alpha}\in D\cap \overline{D_{\zeta}}^{{\varepsilon}}$.
 Then $p_{\alpha}\in \overline{D_{\zeta}}^{{{\tau}_{{\omega}_1}}}$, so $D$ is not discrete. 
 Hence, $s(\mc X)={\omega}$.   
 \end{proof}

\section{Problems.}

 By \cite[12.5]{HandBook-vD},
  if $X$ is a regular, feebly compact, first countable space  with  $|X|<\mf b$, then $X$ is countably compact.

  \begin{problem}In the statement above, is it necessary
    to assume that $X$ is $T_3$? What about $T_2$ spaces?
   \end{problem}

  \begin{problem}
 Is there a    regular, feebly compact, but not countably compact,  first countable space  with  $|X|=\mf b$
  in ZFC?\end{problem}

Theorem  \ref{pr:S-space} shows that   the following question arises  naturally.
\begin{problem}
Does the existence of an S-space imply the existence of a first countable \PNC\  space  with countable spread? 
 \end{problem}

 In Theorem \ref{tm:Thm1.2} we obtain only consistency.
 \begin{problem}
  Is there,  in ZFC,  a 0-dimensional   \Rhered ``\DRC\ but $\neg$\DCC''  $ T_2$ space ?
\end{problem}

Concerning  the next problem, 
we  have a consistency result without assuming the first countability (see Theorem \ref{tm:Thm1.2}). 
\begin{problem}
  Is it consistent that 
there exists a first countable, 0-dimensional $ T_2$,    \Rhered ``''\DRC\ but $\neg$\DCC''  space ?
\end{problem}

We know that there exist arbitrarily large \Rhered ``\DCC\ but $\neg$\CC''  (or ``pseudocompact, but $\neg$\SP'') spaces.

\begin{problem}\label{pr:xxx}
Are there arbitrarily large  \Rhered ``\DRC\ but $\neg$\DCC'' (or ``\SP\ but $\neg$\DRC'') spaces?
  \end{problem}

\bibliographystyle{abbrv}

  \end{document}